\documentclass[12pt]{article}
\usepackage{amscd}
\newcommand{\al}{\alpha}
\newcommand{\be}{\beta}
\newcommand{\ep}{\epsilon}

\newcommand{\dbar}{\overline{\partial}}

\newcommand{\dbr}[1]{\partial_{\overline{#1}}}
\newcommand{\ddt}[1]{\frac{\partial #1}{\partial t}}

\newcommand{\FH}{\hat{F}}

\newcommand{\Tr}{\textrm{Tr} \, }

\newcommand{\so}{\mathcal{O}}

\newcommand{\C}{\mathcal{C}}

\newcommand{\oR}{\mathcal{O}_R}
\newcommand{\se}{\mathcal{E}}
\newcommand{\ser}{\mathcal{E}_R}
\newcommand{\serx}{\mathcal{E}_{R,x}}
\newcommand{\sei}{\mathcal{E}^{\infty}}
\newcommand{\seix}{\mathcal{E}^{\infty}_x}
\newcommand{\shf}{\mathcal{F}}
\newcommand{\shfr}{\mathcal{F}_R}
\newcommand{\shfi}{\mathcal{F}^{\infty}}
\newcommand{\shpr}{\mathcal{P_R}}
\newcommand{\ssr}{\mathcal{S_R}}

\newcommand{\rank}{\textrm{rank} \,}
\newcommand{\ov}[1]{\overline{#1}}
\begin{document}
\title{Complex Frobenius}
\newcounter{theor}
\setcounter{theor}{1}
\newtheorem{theorem}{Theorem}[section]
\newtheorem{proposition}{Proposition}[section]
\newtheorem{lemma}{Lemma}[section]
\newtheorem{defn}{Definition}[theor]
\newtheorem{corollary}{Corollary}[section]
\newenvironment{proof}[1][Proof]{\begin{trivlist}
\item[\hskip \labelsep {\bfseries #1}]}{\end{trivlist}}
\renewcommand{\theequation}{\thesection.\arabic{equation}}

\centerline{\bf A COMPLEX FROBENIUS THEOREM,} 
\centerline{\bf MULTIPLIER IDEAL SHEAVES}
\centerline{\bf AND HERMITIAN-EINSTEIN METRICS}
\centerline{\bf ON STABLE BUNDLES}
\bigskip
\centerline{Ben Weinkove}
\centerline{Department of Mathematics, Columbia University}
\centerline{New York, NY 10027}
\centerline{Email: {\it weinkove@math.columbia.edu}}
\bigskip

\setlength\arraycolsep{2pt}
\addtocounter{section}{1}

\bigskip
\noindent
{\bf 1. Introduction}
\bigskip

In \cite{Ko}, J.J. Kohn introduced the notion of `subelliptic multipliers' in
his work on the $\dbar$-Neumann problem on pseudo-convex domains in order to prove
subelliptic estimates on the boundary under certain conditions.  He believed that
his techniques would have wide applications in deriving estimates by algebraic
methods.  Following this line of approach, Nadel \cite{Na} defined the {\em
multiplier ideal sheaf}, and used it to prove the existence of K\"ahler-Einstein
metrics on certain Fano manifolds with symmetries (his proofs were later
simplified in \cite{DeKo}.)  Given a plurisubharmonic function
$\phi$ on  the manifold, the multiplier ideal sheaf of $\phi$ is the sheaf of
germs of holomorphic functions $f$ satisfying
$$\int |f|^2 e^{-\phi} < \infty.$$
The sheaf measures the singularities of $\phi$.  The concept of the multiplier
ideal sheaf, which had already been considered implicitly in the work of Bomberi
\cite{Bo}, Skoda \cite{Sk} and Siu \cite{S1}, is an important one in algebraic
geometry, and has led to some important breakthroughs, such as the proof by Siu
\cite{S3} of the big Matsusaka theorem and effective results on global generation
and very ampleness (e.g. \cite{S4}, \cite{De}.)

In a number of lectures and seminars (\cite{S5} for example), Y-T. Siu has 
advocated that generalized multiplier ideal sheaves
have many applications to problems in partial differential equations and
geometry.  A multiplier ideal sheaf should measure the position and extent of the
failure of an \emph{a priori} estimate.  A key concept is that the
sheaf $\mathcal{I}$ should satisfy a differential inclusion
relation which arises from known estimates.  In Kohn's
case (\cite{Ko}, Proposition 4.7.(G)), the differential relation
takes the form
$$\partial \mathcal{I} \subset \mathcal{M},$$
where $\mathcal{M}$ is Kohn's `multiplier module'.
In this paper, we look at the notion of multiplier ideal sheaves in the context
of the problem of the existence of Hermitian-Einstein metrics on stable bundles
over compact K\"ahler manifolds.   It is our hope that this may provide some
clues as to how to  proceed in other situations.  The problem of existence of
Hermitian-Einstein metrics has been solved over curves by Narasimhan and
Seshadri \cite{NS} (see also \cite{D1} for an alternative proof), over algebraic
surfaces by Donaldson \cite{D2} and over general compact K\"ahler manifolds by
Uhlenbeck and Yau \cite{UY}.  In \cite{UY}, using a continuity method, it is
shown that proving the existence of an  Hermitian-Einstein metric can be reduced
to proving a certain \emph{a priori} estimate.  If this estimate fails to hold
then they show that there exists a {\em destabilizing subsheaf} which contradicts
the assumption that the original bundle was stable. 

We give an alternative construction of the destabilizing subsheaf in the case of
algebraic surfaces, following the ideas of Siu.  Instead of the continuity method,
we use the Donaldson heat flow (see section 3) which is roughly equivalent.  The
crucial estimate is the
$C^0$ estimate for this flow of metrics.  Abusing terminology, we define a
`multiplier ideal sheaf' picking out those directions in which the metric does
not blow up.  We show that, using a bound on the contracted curvature $\Lambda F$
along the flow, the multiplier ideal sheaf
$\mathcal{F}$ satisfies the differential inclusion relation 
$$ \dbar \mathcal{F} \subset \mathcal{F} \otimes T^{0,1}.$$
The key property that the sheaf must satisfy is coherence.  We show this 
using the differential inclusion relation and a complex Frobenius theorem
which is given as follows. 

Let $E$ be a holomorphic vector bundle of rank $r$ over a complex manifold $X$ 
of complex dimension $n$.  Denote the `partial connection' associated to the 
holomorphic structure on $E$ by $\dbar$.  Denote by $\so$, $\oR$ and $\C^{\infty}$
the sheaves of holomorphic, complex-valued real analytic and smooth
functions on $X$ respectively.  Let\
$\se$, $\ser$ and $\sei$ denote the sheaves of holomorphic, real analytic and
$C^{\infty}$ sections of $E$ respectively. 
$T^{0,1}$ will denote the corresponding sheaf of (0,1) forms on $X$.

\bigskip
\pagebreak[3]
\noindent
{\bf Theorem 1}
{\it
Let $\shf_R \subseteq \ser$ be a
subsheaf satisfying
\begin{enumerate}
\item[(i)] $\shf_R$ is locally finitely generated;
\item[(ii)] $\dbar \shf_R \subseteq \shf_R \otimes T^{0,1}.$
\end{enumerate}
Then $\shf_R$ is locally finitely generated by holomorphic sections of $E$.  Denote
by
$\shf$ the sheaf of $\so$-modules $\ker \dbar|_{\shf_R}$.  It follows that
$\shf \subset \se$ is coherent and
\begin{equation} \label{sheafeqn}
\shf \otimes \oR = \shf_R.
\end{equation}
The corresponding result holds for a subsheaf $\shf^{\infty} \subseteq \sei$ if
condition (i) is replaced by
\begin{enumerate}
\item[(i)'] For each $x$ in $X$, there exists an
open set $U$ containing $x$
 and a free resolution
$$\cdots \longrightarrow (\C^{\infty}_U)^{d_2}
\longrightarrow (\C^{\infty}_U)^{d_1}
\longrightarrow (\C^{\infty}_U)^{d_0} \longrightarrow \shf^{\infty}_U
\longrightarrow 0,$$
for $d_i \ge 0$.
\end{enumerate}
}

\bigskip

If the subsheaf is locally free, this follows from the
Koszul-Malgrange integrability theorem for vector bundles \cite{KM}, which is a
linear version of the Newlander-Nirenberg theorem for integrable almost complex
structures
\cite{NN}.  

A result similar to the $C^{\infty}$ case of Theorem 1 has been proved by very
different methods in a recent preprint of Pali \cite{Pa}.  His theorem applies to
sheaves of smooth functions with a $\dbar$ operator satisfying $\dbar^2=0$, and
which  admit local free resolutions of finite length.  Pali does not assume that the
sheaf is a subsheaf of a holomorphic vector bundle.  In the case that it is,
condition
\emph{(i)'} of Theorem 1 is slightly weaker than his condition, since we do not
assume that the resolution is of finite length. 

Note that the differential inclusion relation \emph{(ii)} is not
enough by itself to prove coherence.  This can be seen by the following
counterexample.  Fix an open set
$U$ in $X$ and consider the subsheaf with stalk $\serx$ ($\seix$) for $x$ in
$U$ and the zero stalk otherwise.  This satisfies condition
\emph{(ii)} but is clearly not coherent.

In our application of Theorem 1, we
will only use the real analytic case.  We define our multiplier ideal sheaf as
follows.
Suppose that $X$ is now a compact K\"ahler manifold, and let $\tilde{X}$ 
be an open subset of $X$.  Let $\{ H_k
\}_{k=0}^{\infty}$ be a sequence of metrics on
$E$ over $\tilde{X}$, and define endomorphisms $h_k = H_0^{-1} H_k$.  Define a
presheaf $\shpr \subseteq
\ser$ over $\tilde{X}$ by setting the sections over an open set 
$U \subset \tilde{X}$ to be
$$\shpr (U) = \{ s \in \ser(U) \ | \ \int_U |s|^2_{H_k} \rightarrow 0 \textrm{ 
as } k \rightarrow \infty \},$$
where integration is performed with respect to the volume form induced from the 
K\"ahler metric on $X$.  Then let $\shf_R$ be the sheaf associated to the presheaf
$\shpr$.   We will call $\shf_R$ the \emph{multiplier ideal sheaf} associated to
the sequence of metrics $\{ H_k \}$.  (For a different definition of a multiplier ideal sheaf for bundles, see
\cite{DC}.)  If the sequence of metrics satisfies
certain conditions, then we can prove the differential inclusion relation and
apply Theorem 1 to obtain a coherent subsheaf of $E$. 

\bigskip
\pagebreak[3]
\noindent
{\bf Theorem 2} 
{\it Suppose that there exists a positive constant
$C$ and a real analytic endomorphism $h_{\infty}$ of $E$ such that 
the sequence of metrics $\{ H_k \}_{k=0}^{\infty}$ satisfies
\begin{enumerate}
\item[(i)] $H_k \le H_0$ for all $k$;
\item[(ii)] $i \hat{F}_{H_k} H_0 \le C H_0$ for all $k$;
\item[(iii)] $h_k \rightarrow h_{\infty}$ uniformly on compact subsets of
$\tilde{X}$, for $h_k = H_0^{-1} H_k.$
\end{enumerate}
Then the multiplier ideal sheaf $\shf_R$ associated to $\{ H_k \}$ defines a
holomorphic coherent subsheaf $\shf \subset \se$ over
$\tilde{X}$.}

\bigskip

We will see how to use Theorem 2 to construct the destabilizing
subsheaf in a proof of the theorem of Donaldson \cite{D2}:

\bigskip
\noindent
{\bf Theorem 3}
{\it Let $X \subset \mathbf{CP}^N$ be a projective algebraic surface and
let $\omega$ be the K\"ahler metric on $X$ induced from the Fubini-Study
metric on $\mathbf{CP}^n$.  Let $E$ be a holomorphic vector bundle of
rank $r$ over $X$ which is stable with respect to the metric $\omega$.
Then $E$ admits a Hermitian-Einstein metric.}

\bigskip
 
It should be emphasized that the proof of Theorem 3 given here differs from
existing proofs only in the construction of the destabilizing subsheaf and does
not provide a simplification (in fact the given proof may be easily shortened.)
However, the method should help elucidate a certain point of view and provide a
useful model for the further study of multiplier ideal sheaves. We
use the Yang-Mills heat flow and the associated flow of metrics, taking estimates
from
\cite{D2}, \cite{DoKr} and \cite{UY}.  The heat flow method has also been used by
Simpson \cite{Si} and de Bartolomeis and Tian \cite{DT} to prove generalizations
of the Theorem of Donaldson-Uhlenbeck-Yau.

The paper is arranged as follows.  In section 2 we prove the complex Frobenius
theorem. In section 3 we prove Theorem 2 by making use of the differential
inclusion relation.  In section 4 we give a proof of Theorem 3, using
the Yang-Mills flow and Theorem 2.

\bigskip
\noindent
{\it Notation}
\bigskip

We will use the notation from \cite{D2}.  In particular, the K\"ahler
form
$\omega$ on $X$ can be written in normal
coordinates at a point as
$$ \omega = \frac{1}{2} i \sum_{j=1}^2 dz^j \wedge d\overline{z}^j,$$
and the operator $\Lambda$ acts on $(1,1)$ forms at this point by
$$\Lambda (\sum_{j,k =1}^2 a_{j\overline{k}} dz^j \wedge d\overline{z}^k) = -2i
\sum_{j =1}^2 a_{j \overline{j}}.$$
If $\phi$ and $\psi$ are form-valued sections of $E$, define
$$\langle \phi , \psi \rangle_H = \sum_{\al, \be =1}^{r} H_{\al \overline{\be}}
\phi^{\al}
\wedge
\overline{\psi^{\be}}.$$ 
Note that if $H_0$ and $H$ are metrics on $E$, and $h= H_0^{-1} H$, then the
associated curvatures $F_{H_0}$ and $F_H$ differ by
$$F_{H} = F_{H_0} + \dbar (h^{-1} \partial_0 h).$$
For convenience, write $\hat{F}$ for $\Lambda F$.

\pagebreak[3]
\setcounter{equation}{0}
\addtocounter{section}{1}
\bigskip
\noindent
{\bf 2. Complex Frobenius theorem}
\bigskip

In this section we give a proof of Theorem 1.  We begin with the real analytic
case.  The statement is local so it will suffice to prove the following. Let
$\shf_R$ be a subsheaf of 
$\oR^r$ over 
the polydisc $U_{\eta} = \{ z \in \mathbf{C}^n \ | \ |z_j| < \eta , 
j=1,\ldots, n \}$ for some $\eta >0$ satisfying the hypotheses $(i)$ and $(ii)$ of
Theorem 1.  Let 
$f_1, \ldots, f_q$ be sections of $\shfr$ over $U_{\eta}
$ generating the stalks $\shf_{R,x}$ for all $x$ in $U_{\eta}$.  We will show that 
there exists some $\eta'$ with $0 < \eta' < \eta$ and sections $g_1, \ldots, 
g_q$ of $\shf_R$ over $U_{\eta'}$ such that $\dbar g_i = 0$ for $i=1, \ldots, 
q$, with the $g_i$ generating $\shf_{R,x}$ for all $x$ in $U_{\eta'}$.

We will give a proof by induction on the number of variables $z_1, \ldots, 
z_p$ in which the generators are holomorphic.  

\bigskip
\noindent
{\it Step 1. }
We begin by first showing that for small enough $\eta'$
there exist sections $g_1, \ldots, g_q$ over $U_{\eta'}$ satisfying $\dbr{1}g_i 
=0$, 
and such that $g_1, \ldots, g_q$ generate $\shf_{R,x}$ for all $x$ in $U_{\eta'}$.  
By 
hypothesis, there exists a $q \times q$ matrix of real analytic functions $A_1$
on $U_{\eta}$  satisfying
$$\pmatrix{\dbr{1} f_1 \cr \vdots \cr \dbr{1} f_q} = A_1  \pmatrix{f_1 \cr 
\vdots \cr f_q}.$$
We will show that for small enough $\eta'$ there exists a 
real analytic map $B: U_{\eta'} \rightarrow GL(q, \mathbf{C})$ satisfying
\begin{equation} \label{eqndiff1}
\dbr{1} B + BA =0.
\end{equation}
To see that this will give the desired conclusion, set
$$ \pmatrix{g_1 \cr \vdots \cr g_q} = B \pmatrix{f_1 \cr \vdots \cr f_q}.$$
Since $B$ is invertible, the $g_1, \ldots, g_q$ generate $\shf_{R,x}$ for $x$ in 
$U_{\eta'}$ and differentiating the above gives
$$ \pmatrix{\dbr{1} g_1 \cr \vdots \cr \dbr{1} g_q} = (\dbr{1} B) \pmatrix{f_1 
\cr \vdots \cr f_q} + B A_1 \pmatrix{f_1 \cr \vdots \cr f_q} = 0.$$
We will pull back the $f_i$ via the dilation map $\delta_r :
\mathbf{C}^n \rightarrow  \mathbf{C}^n$  given by $(z_1, z_2, \ldots, z_n) 
\mapsto (r z_1, z_2, \ldots, z_n)$ for  some 
$r<1$.  Notice that $A_1$ scales by the factor $r$.  We will solve the problem
for these new $f_i$ and then transform back. 

Denote by $D_{\eta}$ the open disc of
radius $\eta$ in $\mathbf{C}$. Multiply $A_1$
by a smooth cut-off function $\psi=\psi(|z_1|)$ with compact support in
$D_{\eta/r}$ and equal to 1 on $D_{\eta/2r}$.
We can then regard $A_1$ as a function on $\mathbf{C}$ with
parameters $z_2, \ov{z}_2, \ldots, z_n, \ov{z}_n$ varying
analytically in a small polycylinder which may be shrunk if necessary.  Replacing
these variables by $\zeta_1,
\ldots, \zeta_{2n-2}$, we can view $A_1$ as having holomorphic
parameters $\zeta_i$.  $A_1$ is smooth on $\mathbf{C}$, has compact support in
$D_{\eta/r}$ and is real analytic in $D_{\eta/2r}$.  By choosing $r$ small enough,
we can suppose that, for a fixed $0 < \ep < 1$, the norm $\| A_1 \|_{C^{\ep}}$ is as
small as we like. 

Now set $B=I+F$, so that equation (\ref{eqndiff1}) becomes
\begin{equation} \label{eqndiff2}
\dbr{1} F + (I+F)A_1 = 0.
\end{equation}
We will use the following implicit function theorem to solve
(\ref{eqndiff2})  for small real analytic $F$.  

\pagebreak[3]
\begin{theorem}
Let $Y_1$, $Y_2$ and $Z$ be Banach spaces, and set $Y=Y_1 \times Y_2$.  Let 
$\Phi:Y \rightarrow Z$ be a smooth (holomorphic) map.  Fix a point $(y_1, y_2)$ in
$Y$.  If  the partial derivative $(D_2 \Phi)(y_1,y_2) : Y_2 \rightarrow Z$ is
surjective  and admits a bounded right inverse $P: Z \rightarrow Y_2$, then there is
a  smooth (holomorphic) map $f$ from a neighbourhood of $y_1$ in $Y_1$ to a
neighbourhood of
$y_2$ in $Y_2$ such that
$$ \Phi( \eta, f(\eta)) = \Phi (y_1, y_2).$$ 
It follows that if $\eta_{z}$ is a family of elements of $Y_1$ in the above
neighbourhood of $y_1$, varying smoothly (holomorphically) in a parameter $z$, then
$f(\eta_{z})$ will also vary smoothly (holomorphically) in $z$.
\end{theorem}
          
Apply the theorem with $Y_1$ and $Z$ the Banach spaces of 
$q \times q$ matrix-valued functions on $\mathbf{C}$ with the
norm
$\|\
\|_ {C^\ep}$ bounded, and $Y_2$ the same space with norm $\| \
\|_ {C^{1 + 
\epsilon}}$.  Take $(y_1,y_2)=(0,0)$ and 
$$ \Phi(A_1,F) = \dbr{1}F + (I+F)A_1.$$
Then $(D_2 \Phi)(y_1,y_2): Y_2 \rightarrow Z$ is given by
$$(D_2 \Phi) (y_1, y_2) (\al) = \dbr{1}\al,$$
which has bounded right inverse $P: Z \rightarrow Y_2$ given by
$$P(\beta)(z) = \frac{1}{2\pi i} \int_{\mathbf{C}} \frac{\be(w)}{w-z} 
dw\wedge d\overline{w}.$$
Hence, if $\| A_1 \|_{C^{\epsilon}}$ is small enough then there exists $F=F(A_1)$ 
with $\| F \|_{C^{1+ \epsilon}}$ small satisfying equation (\ref{eqndiff2}).
Now since $A_1$ depends holomorphically on the parameters $\zeta_i$, and since
$\Phi$ is holomorphic, the solution $F$ will also be holomorphic in these parameters.
By making the 
norm of $A_1$ smaller if necessary, the solution $B=I+F$ will be invertible.   
Moreover, since $F$ is a solution of the elliptic equation
$$ (\partial_{\overline{z}_1} + \partial_{\overline{\zeta}_1} + \ldots + 
\partial_{\overline{\zeta}_{2n-2}})F + (I+F)A_1 =0,$$
it follows that $F$ must also be  
analytic in the variables $z_1$ and $\overline{z}_1$ by the theorem that 
sufficiently smooth solutions to elliptic equations with real analytic
coefficients are themselves real analytic.
Replacing the $\zeta_1, \ldots, \zeta_{2n-2}$ with $z_2,\ldots, \ov{z}_n$, we see
that $F$ is analytic in $z_1, \ov{z}_1, \ldots, z_n, \ov{z}_n.$  We 
have now completed the first 
step in the induction.

\bigskip
\noindent
{\it Step 2. }  Fix $p$ with $1 \le p \le n-1$.  Our inductive hypothesis is 
that for $\eta'$ small enough, there exist sections $f_1, \ldots, f_q$ over $U_
{\eta'}$ satisfying $\dbr{j} f_i=0$ for $j=1, \ldots, p$ and generating 
$\shf_{R,x}$ for $x$ in $U_{\eta'}$.  We will show that, for some $\eta''>0$,  there
exist sections $g_1, \ldots, g_q$ over $U_{\eta''}$ satisfying 
$\dbr{j} g_i =0$ for $j= 1, \ldots, p+1$ and such that $g_1, \ldots, g_q$ 
generate $\shf_{R,x}$ for $x$ in $U_{\eta''}$.

We know that there exist real analytic functions of $q \times q$ matrices $A_k$ 
on $U_{\eta'}$
for $k=1,\ldots, n$ satisfying
\begin{equation} \label{eqndbar}
\pmatrix{\dbr{k} f_1 \cr \vdots \cr \dbr{k} f_q} = A_k \pmatrix{f_1 \cr \vdots 
\cr f_q}.
\end{equation}
The inductive hypothesis implies that $A_k=0$ for $k=1, \ldots, p$.  We will 
show that, for $\eta''$ small enough, there exists a real analytic function 
$B:U_{\eta''} \rightarrow GL(q, \mathbf{C})$, holomorphic in the variables $z_1, 
\ldots, 
z_p$, satisfying 
$$\dbr{p+1}B + B A_{p+1} = 0.$$
Then, as before, the sections $g_i$ given by
$$ \pmatrix{g_1 \cr \vdots \cr g_q} = B \pmatrix{f_1 \cr \vdots \cr f_q}$$
will satisfy the required properties.

Before using the implicit function theorem, we must ensure that $A_{p+1}$ is 
holomorphic in $z_1, \ldots, z_p$.  Notice that in equation (\ref{eqndbar}) for
$k=p+1$,
$A_{p+1}$ is an analytic function of all of the variables $z_1, \overline{z}_1,
\ldots, z_n,
\overline{z}_n$ whereas the $\dbr{p+1} f_i$ and the $f_i$ do not depend on
$\overline{z}_1, \ldots, \overline{z}_p$.  Setting
$\overline{z}_1 = \ldots = \overline{z}_p = 0$ we get an $A_{p+1}$ which satisfies
(\ref{eqndbar}) and is holomorphic in $z_1, \ldots, z_p$.

We can now finish the proof.  Let $B = I+F$.  We use the same method as before 
to solve
$$\dbr{p+1} F + (I+F)A_{p+1} = 0$$
for $F$ small, real analytic, and holomorphic in $z_1, \ldots, z_p$, given that 
$A_{p+1}$ is 
holomorphic in $z_1, \ldots, z_p$.  
After dilating the $z_{p+1}$ coordinate and multiplying by a cut-off function 
$\psi=\psi(|z_{p+1}|)$, we 
can take $A_{p+1}$ to be a smooth function on $\mathbf{C}$ in the $z_{p+1}$ 
variable, with the $z_i$ ($i \neq {p+1}$) regarded as parameters which vary in 
some small polycylinder.  As before we can assume that $\| A_{p+1} 
\|_{C^{\epsilon}}$ is as small as we like. Apply 
the implicit function theorem to 
get a solution $F$ with $\| F \|_{C^{1+\epsilon}}$ small and which is  
analytic in $z_{p+1}, \overline{z}_{p+1}, \ldots, z_n, \overline{z}_n$ and
holomorphic in the variables $z_1, \ldots, z_p$.  This completes the induction, and
the theorem for the real analytic case follows.

\bigskip

For the smooth case, we cannot apply the argument used in Step
2 to get an $A_{p+1}$ holomorphic in $z_1, \ldots, z_p$.  We
will replace this with another induction argument.  Define a 
proposition $P_k$ for $1 \le k \le p$ as follows.

\bigskip
\noindent
{\bf Proposition $\mathbf{P_k}$ } \emph{Let
$\mathcal{G^{\infty}} \subset (\C^{\infty})^d$ be a subsheaf over $U_{\delta}$
for some $\delta>0$, satisfying condition (i)' of the theorem.  Suppose that
$\mathcal{G^{\infty}}$ is generated over $U_{\delta}$ by sections
$r_1,\ldots, r_l$, which are holomorphic in $z_1, 
\ldots, z_p$.  Let $a_1, \ldots, a_l$ be smooth functions
on
$U_ {\delta}$ 
satisfying
$$ \pmatrix{\dbr{j} a_1 \cdots \dbr{j}a_l} \pmatrix{r_1 \cr \vdots \cr r_l} = 0, 
\qquad \textrm{for } j=1, \ldots, k.$$
Then
there exists $\delta'$ with $0<\delta' < \delta$ and $\pmatrix{s_1 \cdots s_l}$ 
in the $\C^{\infty}$ sheaf of relations $\mathcal{R}(r_1, \ldots, r_l)$ over
$U_{\delta'}$ with 
$$ \dbr{j}(a_1 - s_1) = \ldots = \dbr{j}(a_l - s_l) = 0, \qquad \textrm{for } 
j=1, \ldots, k.$$}

We will show that $P_k$ holds for $k=p$.  Then since
$$ \dbr{j} A_{p+1} \pmatrix{f_1 \cr \vdots \cr f_q} = 0, \qquad \textrm{for } j=1,
\ldots, p,$$
we can apply this result to each row of $A_{p+1}$.  Hence, after
making permissible changes to $A_{p+1}$ and 
passing to a smaller open set, we can assume that $A_{p+1}$ is holomorphic in $z_1,
\ldots, z_p$.   We will need to use the main inductive hypothesis in the proof of
this induction.

\bigskip
\noindent
{\bf Proof of $\mathbf{P_1}$ } By assumption, the sheaf of
relations
$\mathcal{R} (r_1,\ldots, r_l)$ is locally finitely generated. Suppose the
generators are 
$$ h_1 = \pmatrix{h_{11} \, h_{12} \cdots h_{1l}}, \ldots, h_{m}=\pmatrix{h_{m 1}
\, h_{m 2} 
\cdots h_{m l}}.$$
Since $r_1, \ldots, r_l$ are holomorphic in $z_1$, and $\mathcal{R} (r_1, \ldots,
r_l)$ satisfies
\emph{(i)'}, we can apply the
\emph{main} inductive hypothesis in the case $p=1$ to this subsheaf.  Hence we can
assume that, after choosing 
$\delta'$ small enough, the 
generators $h_1, \ldots, h_m$
are holomorphic in the variable $z_1$.  Now there exist smooth
functions 
$t_1, \ldots, t_m$ on $U_{\delta'}$ satisfying
$$\pmatrix{\dbr{1} a_1 \cdots \dbr{1} a_l} = t_1 h_1 
+ \ldots + t_m h_m.$$
Let $\hat{t}_1, \ldots, \hat{t}_m$ be smooth functions
satisfying
$\dbr {1}\hat{t}_i = t_i$ (shrinking $\delta'$ slightly.) 
Define
$$\pmatrix{s_1 \cdots s_l} = \hat{t}_1 h_1 + \ldots 
+ \hat{t}_m h_m.$$
Then we see that $\pmatrix{s_1 \cdots s_l} \in \mathcal{R}(r_1, \ldots, r_l)$ 
and 
$$ \dbr{1}(a_1 - s_1) = \ldots = \dbr{1}(a_l - s_l) = 0,$$
completing the proof of $P_1$.

\bigskip
\noindent
{\bf Proof of $\mathbf{P_k \Rightarrow  P_{k+1}}$ }  We assume $1 \le k \le
p-1$.  Suppose that
$$ \pmatrix{\dbr{j} a_1 \cdots \dbr{j}a_l} \pmatrix{r_1 \cr \vdots \cr r_l} = 0, 
\qquad \textrm{for } j=1, \ldots, k+1.$$
By the inductive hypothesis, we can suppose, after making some allowed 
changes and taking $\delta'$ small enough, that $\dbr{j}a_i=0$ for $j=1,\ldots,
k$  and 
$i=1,\ldots, l$.  As before let the generators of $\mathcal{R}(r_1,\ldots, r_l)
$ be
$$ h_1 = \pmatrix{h_{11} \, h_{12} \cdots h_{1l}}, \ldots, h_m = \pmatrix{h_{m1} \,
h_{m2} 
\cdots h_{ml}}.$$
Since $r_1, \ldots, r_l$ are holomorphic in $z_1, \ldots, z_k$, by applying 
the main inductive hypothesis, we can assume that $h_1, \ldots, h_m$ are 
holomorphic in $z_1, \ldots, z_k$.  There exist $t_1, \ldots, t_m$
smooth functions with
$$\pmatrix{\dbr{k+1} a_1 \cdots \dbr{k+1}a_l} = t_1 h_1 + \ldots + t_m h_m.$$
Notice that by applying $\dbr{j}$ we get that
$$\pmatrix{\dbr{j} t_1 \cdots \dbr{j} t_m} \pmatrix{h_1 \cr \vdots \cr h_m} = 0, 
\qquad \textrm{for } j=1, \ldots, k.$$
We can apply the \emph{current} inductive hypothesis to see that, after shrinking 
$\delta'$, there exists $\pmatrix{s'_1 \cdots s'_m}$ 
in $\mathcal{R}(h_1, \ldots, h_m)$ such that
$$\dbr{j}(t_1 - s'_1) = \ldots = \dbr{j}(t_m - s'_m) = 0, \qquad \textrm{for } 
j=1, \ldots, k.$$
By replacing the $t_i$ by $t_i - s'_i$, we can then assume that $\dbr{j} t_i =0
$ 
for $j=1, \ldots , k$.  Let $\hat{t}_i$ solve $\dbr{k+1} \hat{t}_i = t_i$ 
for 
$i=1, \ldots, m$ where the $\hat{t}_i$ are smooth and
holomorphic in
$z_1, 
\ldots, z_k$.  Then, as before, set
$$\pmatrix{s_1 \cdots s_l} = \hat{t}_1 h_1 + \ldots + \hat{t}_m h_m.$$
Then $\pmatrix{s_1 \cdots s_l} \in \mathcal{R}(r_1 ,\ldots,
r_l)$ and
$$\dbr{j}(a_1 - s_1) = \ldots = \dbr{j}(a_l - s_l)=0, \qquad \textrm{for } 
j=1,\ldots, k+1.$$
This proves $P_{k+1}$.

It is now straightforward to complete the proof in
the smooth case.  To see that $\shf = \ker \dbar|_{\shfi}$ is coherent and to
establish
$$ \shf \otimes \C^{\infty} = \shfi,$$
one can use the result of Malgrange
\cite{Ma} on the flatness of $\C^{\infty}$ over $\so$, or apply Corollary 6.3.6 of
\cite{Ho}.

\pagebreak[3]
\setcounter{lemma}{0}
\setcounter{theorem}{0}
\setcounter{equation}{0}
\addtocounter{section}{1}
\bigskip
\bigskip
\noindent
{\bf 3. The multiplier ideal sheaf} 
\bigskip

In this section, we give a proof of Theorem 2.  Let
$H_k$ be a sequence of metrics as given in the hypothesis of the theorem,
defining the multiplier ideal sheaf $\shf_R$ over $\tilde{X}$.  We now prove the
differential inclusion relation. 

\addtocounter{lemma}{1}
\begin{proposition} \label{propdbar}
$\dbar \shf_R \subseteq \shf_R \otimes T^{0,1}$.
\end{proposition}
\begin{proof}
Let $x$ be in $\tilde{X}$, and let $U$ be an open set in $\tilde{X}$
containing $x$.  Choose 
$V$ compact and $W$ open such that $x \in W \subset V \subset U$ and let
$\psi$  be a smooth cut-off function supported in $V$ and equal to 1 in
$W$.  Suppose that $s$ satisfies
$$ \int_U |s|^2_{H_k} \rightarrow 0 \quad \textrm{as } k \rightarrow \infty.$$
Then
\begin{eqnarray*}
\int_W |\dbar s |^2_{H_k} & = & 
-i \int_W \Lambda \langle \dbar s, \dbar s \rangle_{H_k} \\
& \le & -i \int_X \Lambda
\langle \psi h_k \dbar s, \dbar s \rangle_{H_0} \\
& = & -i \int_X \Lambda \langle \psi (\partial_0 h_k) \dbar s, s
\rangle_{H_0} - i \int_X \Lambda 
\langle \psi h_k \partial_0 \dbar s, s \rangle_{H_0} \\
&& \mbox{} - i \int_X \Lambda \langle (\partial \psi) h_k \dbar s, s
\rangle_{H_0} \\ & = & -i \int_X \Lambda \langle \psi (h_k)^{-\frac{1}{2}}
(\partial_0 h_k) \dbar s,  (h_k)^{\frac{1}{2}} s \rangle_{H_0} \\
&& \mbox{} - i \int_X \Lambda \langle \psi (h_k)^{\frac{1}{2}} \partial_0
\dbar s, (h_k) ^{\frac{1}{2}} s\rangle_{H_0} \\
&& \mbox{} -i \int_X \Lambda \langle (\partial \psi) (h_k)^{\frac{1}{2}}
\dbar s, (h_k) ^{\frac{1}{2}} s\rangle_{H_0} \\
& \le & C (\int_V |(h_k)^{-\frac{1}{2}} \partial_0 h_k |^2_{H_0})^{\frac{1}
{2}} (\int_U |s|^2_{H_k})^{\frac{1}{2}} \\
&& \mbox{} + (\int_V | \partial_0 \dbar s|^2_{H_k} )^{\frac{1}{2}} (\int_U 
|s|^2_{H_k})^{\frac{1}{2}} \\
&& \mbox{} + C (\int_V |\dbar s |^2_{H_k})^{\frac{1}{2}} (\int_U |s|^2_{H_k})^
{\frac{1}{2}} 
\end{eqnarray*}
where $C$ denotes a constant that does not depend on $k$.  The last two terms 
in the last line clearly tend to zero as $k$ tends to zero.  To see that the 
first term in the last line also tends to zero, observe that
\begin{eqnarray*}
\int_X |(h_k)^{-\frac{1}{2}} \partial_0 h_k|^2_{H_0} & = & i \int_X
\Lambda \langle  (h_k)^{-1} \partial_0 h_k, \partial_0 h_k \rangle_{H_0} \\
& = &  i \int_X \Lambda \langle \dbar ((h_k)^{-1} \partial_0 h_k) , h_k
\rangle_{H_0}
\\ & = & \int_X  \langle (i \hat{F}_{H_k} - i \hat{F}_{H_0}), h_k \rangle_{H_0}
\\ & \le & C
\end{eqnarray*}
where we have used the fact that the contracted curvature $i \hat{F}_{H_k}$ is
bounded from above.  This completes the proof.
\end{proof}

In order to apply Theorem 1 we must see that $\shf_R$ is
locally finitely generated.  This follows from the next lemma, which gives an
alternative description of the multiplier ideal sheaf.

\begin{lemma} \label{slemma} Define a sheaf $\ssr$ on $\tilde{X}$ by
$$\ssr (U) = \{s\in \ser(U) \ | \ h_{\infty} s =0 \}.$$
Then $\shf_R = \ssr$.
\end{lemma}
\begin{proof}
Let $x$ be in $\tilde{X}$, and let $x \in U \subset \tilde{X}$.  
Choose $V$ compact and
$W$ open such that $x \in W 
\subset V \subset U$.  If $s$ is in $\ssr(U)$ then
\begin{eqnarray*}
\int_W | s|^2_{H_k} & = & \int_W \langle (h_k - h_{\infty})s, s \rangle_{H_0} 
\\
& \le & \| h_k - h_{\infty} \|_{C^0(V)} \int_V |s|^2_{H_0} \\
& \rightarrow & 0,
\end{eqnarray*}
as $k$ tends to infinity.  This proves $\ssr \subset \shfr$.  On the other hand,
suppose that
$$\int_U |s|^2_{H_k} \rightarrow 0, \qquad \textrm{as } k \rightarrow \infty.$$
Then for any smooth section $t$ of $E$ over $U$,
\begin{eqnarray*}
\int_V \langle h_{\infty} s, t \rangle_{H_0} & = & \int_V \langle (h_{\infty} - 
h'_k) s, t \rangle_{H_0} + \int_V \langle s, t\rangle_{H_k} \\
& \le & \|h_{\infty} - h'_k \|_{C^0(V)} (\int_V
|s|^2_{H_0})^{\frac{1}{2}}(\int_V |t|^2_{H_0})^{\frac{1}{2}}  
 + (\int_U |s|^2_{H_k})^{\frac{1}{2}} (\int_V |t|^2_{H_k})^{\frac{1}{2}} \\
& \rightarrow & 0
\end{eqnarray*}
as $k$ tends to infinity.  Hence $h_{\infty}s = 0$ on $W$ and the lemma is
proved.  
\end{proof}

Since $h_{\infty}$ is real analytic, $\shfr$ is locally the kernel of a
homomorphism between coherent sheaves of $\oR$ modules, and so $\shfr$ is locally
finitely generated.  Applying Theorem 1 completes the proof.

\setcounter{theorem}{0}
\setcounter{lemma}{0}
\setcounter{equation}{0}
\pagebreak[3]
\addtocounter{section}{1}
\bigskip
\noindent
{\bf 4. Hermitian Einstein metrics on stable bundles}
\bigskip

For this section, we will assume the hypotheses of Theorem 3. 
The background 
material for this section can be found in \cite{D2} and \cite{DoKr} (see also
\cite{S2} for a good exposition.)  

Let $H_0$ be a \emph{real analytic}
Hermitian metric on the fibres of $E$.  We can obtain such a metric in
the following way.  Tensor $E$ with the hyperplane line bundle $L$
raised to a high power
$K$ so that the holomorphic sections of $E \otimes L^K$ embed
$X$ into a Grassmannian $G(r,M)$ of $r$-planes in $\mathbf{C}^M$, where $M$ is
the dimension of $H^0(X, E\otimes L^K)$.  Then, via this embedding, pull back the
canonical metric on the universal bundle
$U(r)$ over the Grassmannian to get a real analytic metric on $E
\otimes L^K$.  Then use the canonical real
analytic metric on $L$ to get a real analytic metric $H_0$ on $E$.

Let $\lambda$ be the average of $(\Tr \hat{F}_0)/r$.  We consider the 
flow of endomorphisms $h_t$ with $h_0 = Id_E$ given by the flow 
$$ \ddt{h_t} = -2i h_t (\hat{F}_t - \lambda I)$$
introduced by Donaldson \cite{D2} ($F_t$ is the curvature of the metric $H_t =
H_0 h_t$.)  We call this the Donaldson heat flow.  
It is shown in \cite{D2} that solutions
exist for all time.   The flow is, roughly  speaking, the same as the Yang-Mills
flow up to gauge.  In particular, if $A_0$  is the unitary connection associated
to $H_0$ and the holomorphic  structure on $E$, and if $g_t$ in the
complexified gauge group
$\mathcal{G}^C$  satisfies $g_t^* g_t = h_t$ then $A_t = g_t(A_0)$ satisfies the
gauge  equivalent Yang-Mills flow
$$\ddt{A_t} = -d_{A_t}^* F_{A_t} + d_{A_t} (\al(t))$$
for $\al(t) \in \Omega^0 (\mathbf{g}_E)$. 

Make a 
(real analytic) conformal change to $H_0$ so that $\det h_t =1$ along the flow.  
Also, normalize the metric on $X$ so that $\textrm{Vol}\, (X) =2\pi.$

Suppose first that $\sup_X |h_t|_{H_0}$ is bounded uniformly in $t$.  By
the argument in \cite{D2}, we have a $C^1$ bound and $L^p_2$ bounds for all
$p$ for $h_t$.  Hence, for a sequence of times $t_k \rightarrow \infty$,
$H_k=H_{t_k}$ converges strongly in
$L^p_1$ to a Hermitian metric $H_{\infty}$, which is itself in $L^p_2$.  Then
$$\hat{F}_{H_k} - \lambda I \rightharpoonup \hat{F}_{H_{\infty}} - \lambda I$$
weakly in $L^p$.  Now, Donaldson's functional 
(see \cite{D2}) is bounded below when $E$ is stable (in fact, semi-stable) and 
calculating the functional along the flow gives
$$\hat{F}_{H_k} - \lambda I \rightarrow 0$$
strongly in $L^2$.  Hence we get that $\hat{F}_{\infty} = \lambda I$ 
and by elliptic estimates, $H_{\infty}$ is a smooth (and in
fact real analytic) Hermitian-Einstein metric.  

Suppose then that the $C^0$ norm of $h_t$ is unbounded along the flow.  We will 
show that this leads to a contradiction.  Take a sequence $h_{t_k}$ with the
$C^0$  norm tending to infinity as $k$ tends to infinity.  Then define normalized
endomorphisms $h'_k$ by
$$h'_k = \frac{h_{t_k}}{\sup_X |h_{t_k}|_{H_0}},$$
and let $H'_k = H_0 h'_k$.  
Then we have the 
following theorem, the proof of which is, for the most part, contained in
\cite{D2} and \cite{DoKr}.

\addtocounter{lemma}{1}
\begin{theorem} \label{theoremhinfty}
There exists a finite set of points $\{ x_1, \ldots, x_p \} \subset
X$ such that, after passing to a subsequence, $h'_k$ converges uniformly
on compact sets in $X \backslash \{ x_1, \ldots, x_p \}$ to a real
analytic endomorphism $h_{\infty}$.
\end{theorem}

\begin{proof}
Define $A_k = h'^{1/2}_k (A_0)$.  Then we have
\begin{enumerate}
\item[(i)] $F_{A_k}$ are bounded in $L^2$;
\item[(ii)] $\hat{F}_{A_k}$ are uniformly bounded;
\item[(iii)] $\nabla_{A_k} \hat{F}_{A_k} \rightarrow 0$ in
$L^2$.
\end{enumerate} 
(i) follows from the fact that the Yang-Mills flow is the
gradient flow for the Yang-Mills functional, (ii) follows from
a maximum principle argument (see \cite{D2}) and a proof of (iii) is given in
\cite{DoKr} (Proposition 6.2.14).

It is shown in \cite{DoKr} that there exists a set of points $\{ x_1, \ldots,
x_p \}$ in $X$ and a cover of $X \backslash \{ x_1, \ldots, x_p \}$ by a
system of balls such that, after passing to a subsequence, each connection $A_k$
has curvature with
$L^2$ norm less than a small $\epsilon >0$.  If this $\epsilon$ is small
enough, then by a theorem of Uhlenbeck's \cite{Uh}, we can put the
connections in Coulomb gauge over these balls.  By elliptic estimates and
a standard patching argument (see \cite{DoKr} for details), there exist gauge
transformations $v_k$ over the punctured manifold such that $v_k
h'^{1/2}_k (A_0)$ are bounded in $L^2_2$.  Put $g_k = v_k
h'^{1/2}_k$, and let $\tilde{A}_k = g_k(A_0)$.
Then if $\tilde{A}_k^{1,0}$ and $\tilde{A}_k^{0,1}$
are the $(1,0)$ and
$(0,1)$ components of $\tilde{A}_k$ we have
\begin{eqnarray*}
\tilde{A}_k^{1,0} & = & A_0^{1,0} + g_k^{-1} (\partial_{A_0} g_k) \\
\tilde{A}_k^{0,1} & = & A_0^{0,1} + (\overline{\partial} g_k )g_k^{-1}.
\end{eqnarray*}
Since $g^*_k g_k = h'_k$ we have that the $g_k$ are bounded uniformly,
and from the above, they are bounded in $L^2_3$.  Hence, after passing to a
subsequence, the
$g_k$ converge uniformly on compact subsets to $g_{\infty}$, from which it
follows that $h'_k$ converges uniformly on compact subsets to
$h_{\infty}$.  It remains to see that $h_{\infty}$ is real analytic.  Now
$\tilde{A}_k$ converges over the punctured manifold to $A_{\infty}$, weakly in
$L^2_2$, and by condition (iii) we must have
$$\nabla_{A_{\infty}} \hat{F}_{A_{\infty}} = 0,$$
which is elliptic in a Coulomb gauge.  The elliptic estimates give that
$A_{\infty}$ is smooth.  Since the coefficients of the equation are real
analytic (using the fact that $\omega$ is real analytic), it follows that
$A_{\infty}$ is real analytic.  Now
\begin{eqnarray*}
A_{\infty}^{1,0} & = & A_0^{1,0} + g_{\infty}^{-1} (\partial_{A_0} g_{\infty}) 
\\
A_{\infty}^{0,1} & = & A_0^{0,1} + (\overline{\partial} g_{\infty} )g_{\infty}^
{-1}
\end{eqnarray*}
and $A_0$ is itself real analytic.  It follows that $g_{\infty}$ is real
analytic, and therefore $h_{\infty}$ is also. 
\end{proof}

Now let $\shfr$ be the multiplier ideal sheaf associated to the sequence of
metrics $\{ H'_k \}$ over $\tilde{X} = X \backslash \{x_1, \ldots, x_p \}$. 
Then, since $\hat{F}$ is bounded along the flow, and by the above, we can apply
Theorem 2 to get a coherent sheaf $\shf \subset \se$ over $\tilde{X}$.  In fact,
this defines a coherent subsheaf of $E$ over the whole of $X$.  To see this, 
recall that a coherent subsheaf may be defined in the following way
(see \cite{UY}, section 7, for example.) Let $U_i$ be a cover of $X$ so that
$E|_{U_i}$ is trivial and let
$t_{ij}$ be the transition functions of $E$ defined on the intersections $U_i
\cap U_j$.  Then a coherent subsheaf of $E$ of rank $k$ can be defined to be a set
$\{ f_i \}$ of \emph{rational maps} from $U_i$ to the Grassmannian
$G(k,r)$ with transformation maps $t_{ij}$ on the overlaps.  Thus, in particular,
we can extend over a finite number of points in $X$. 

We need the following lemma.

\begin{lemma}
There exists a positive constant $c$, independent of $k$, such that
$$ \| \emph{Tr} \, h_k \|_{L^1(X)} \ge c \sup_X \emph{Tr} \, h_k.$$
\end{lemma}
\begin{proof}
We use the inequality
$$ \triangle \log \Tr h_k \ge - \frac{1}{2}(|\FH_{H_0}| + |\FH_{H_k}|)$$
from [S2], and the fact that $\hat{F}$ is uniformly bounded along the 
flow. We know that 
the Green's function $G=G(x,y)$ is bounded below by some constant $-A$.
Then 
\begin{eqnarray*}
\log \Tr h_k (x) & = & \frac{1}{2\pi}\int_X \log \Tr h_k (y) dV(y) \\
&& \mbox{} +
\int_X (- \triangle \log \Tr h_k)(y) (G(x,y) + A)dV(y) \\
& \le & \frac{1}{2\pi}\int_X \log \Tr h_k(y) dV(y) + C \\
& \le & \log (\frac{1}{2 \pi} \int_X \Tr h_k(y)dV(y) ) + C,
\end{eqnarray*}
where $C$ is independent of $k$.  This proves the lemma.
\end{proof}

Since $\det h_k =1$ along the flow, the above shows that at least
one eigenvalue of $h_k$ must tend to infinity on a non-empty open set of
$X$.  Hence
$$ 0 < \rank \shf < \rank \se.$$
We will now show that
$\shf$ is destabilizing,  which will give the contradiction.  Note that
$\shf$ is a subbundle outside a subvariety $Z$ of codimension 1.  We will
calculate the second  fundamental form, following the method of \cite{UY}. 
First observe that  given a semi-positive self-adjoint (with respect to
$H_0$) endomorphism
$h$,  one can define an endomorphism
$h^{\sigma}$ by diagonalizing $h$ at a point in  a unitary frame with
respect to $H_0$ and then raising each matrix entry to the power
$\sigma$.  Then note that outside $Z \subset X$, 
the pointwise limit 
$$\pi = \lim_{\sigma \rightarrow 0} (I - h_{\infty}^{\sigma})$$
is the orthogonal projection of $E$ onto the subbundle defined by $\shf$. 
Define the {\em slopes} $\mu(\shf)$ and $\mu(\se)$ by
$$\mu(\se) = \frac{i}{2\pi} \frac{\int_X \Tr(\hat{F}_0)}{\rank(\se)}
\qquad \textrm{and} \qquad
\mu (\shf) = \frac{i}{2 \pi} \frac{\int_X \Tr
(\hat{F}_0(\shf))}{\rank(\shf)}.$$
Note
that 
$\mu(\se) = i \lambda$. 
The following 
lemma completes the proof of Theorem 3.

\begin{lemma} \label{destabilizinglemma}
$$ \mu(\shf) \ge \mu(\se).$$
\end{lemma}
\begin{proof}
First observe that for $0 < \sigma \le 1$,
\begin{eqnarray*}
\int_X |d h'^{\sigma}_k|^2_{H_0} & = & 2 \int_X | \partial_0 h'^{\sigma}_k|^2_
{H_0} \\
& \le & \int_X | (h'_k)^{-\sigma/2} ( \partial_0 h'^{\sigma}_k) |^2_{H_0} \\
& \le & \int_X \langle h_k^{-1} (\partial_0 h_k), \partial_0 h'^{\sigma}_k 
\rangle_{H_0} \\
& = & \int_X i \langle (\hat{F}_k - \hat{F}_0), h'^{\sigma}_k \rangle_{H_0} \\
& \le & C
\end{eqnarray*}
where in the third line, we have used the pointwise inequality
\begin{equation} \label{inequalitysigma}
| h^{-\sigma/2}(\partial_0 h^{\sigma})|_{H_0}^2 \le \langle h^{-1} (\partial_0 
h), \partial_0 h^{\sigma} \rangle_{H_0}
\end{equation}
for $0 < \sigma \le 1$ and $h$ a positive self-adjoint endomorphism (see 
\cite{UY}, Lemma 4.1).
Hence
$$ I - h'^{\sigma}_k \rightharpoonup \tilde{\pi}$$
weakly in $L^2_1$, for some $\tilde{\pi}$ in $L^2_1$, taking the limit as 
$k$ tends to infinity and then as $\sigma$ tends to zero.  Notice that 
$\pi$ and 
$\tilde{\pi}$ must agree on $X \backslash Z$.  
Then
\begin{eqnarray*}
i \int_X \Tr \hat{F}_0 (\shf) & = & i \int_X \Tr ( \hat{F}_0  \pi) -
\int_X  |\partial_0 \pi|_{H_0}^2
\\
& = & i \int_X \Tr ( ( \hat{F}_0 - \lambda I)\pi) + \int_X \Tr (\mu(\se)
\pi) - 
\int_X
|\partial_0 \pi|_{H_0}^2 \\
& = & i \int_X \Tr((\hat{F}_0 - \lambda I)\pi) - \int_X |\partial_0 \pi|_{H_0}
^2 +
\textrm{Vol} \, (X) \, \rank(\shf) \mu(\se).
\end{eqnarray*}
Hence $\mu(\shf) \ge \mu(\se)$ if and only if
$$i \int_X \Tr ((\hat{F}_0- \lambda I)\pi) \ge \int_X |\partial_0 \pi|_{H_0}
^2.$$
First, observe that using the inequality (\ref{inequalitysigma}),
\begin{eqnarray*}
i \int_X \Tr ((\hat{F}_0 - \lambda I) h'^{\sigma}_k) & = & i \int_X \Tr((\hat{F}
_k - \lambda
I) h'^{\sigma}_k) - \int_X \langle \dbar (h_k^{-1}(\partial_0 h_k)), h'^
{\sigma}_k
\rangle_{H_0} \\
& \le & i \int_X \Tr ((\hat{F}_k - \lambda I) h'^{\sigma}_k ) - \int_X 
|\partial_0
h'^{\sigma}_k |_{H_0}^2.
\end{eqnarray*}
Then by lower semicontinuity,
\begin{eqnarray*}
\int_X |\partial_0 \pi|_{H_0}^2 & \le & \liminf_{\sigma \rightarrow 0} \liminf_
{k
\rightarrow \infty} \int_X |\partial_0(I - h'^{\sigma}_k)|_{H_0}^2 \\
& = & \liminf_{\sigma \rightarrow 0} \liminf_{k \rightarrow \infty} \int_X
|\partial_0 h'^{\sigma}_k|_{H_0}^2 \\
& \le & \liminf_{\sigma \rightarrow 0} \liminf_{k \rightarrow \infty} ( -
i \int_X
\Tr ((\hat{F}_0 - \lambda I) h'^{\sigma}_k) + i \int_X \Tr((\hat{F}_k - \lambda 
I)
h'^{\sigma}_k) )\\
& = & i \int_X \Tr ((\hat{F}_0 - \lambda I) \pi) +
\liminf_{\sigma \rightarrow 0} \liminf_{k \rightarrow \infty} i \int_X \Tr
((\hat{F}_k - \lambda I) h'^{\sigma}_k),
\end{eqnarray*}
using the fact that $\int_X \Tr (\hat{F}_0 - \lambda I) = 0$ and $I - h'^
{\sigma}_k
\rightharpoonup \tilde{\pi}$ weakly in $L^2_1$.  Now since the Donaldson 
functional is
bounded below, we know that $\int_X |\hat{F}_k - \lambda I|_{H_k}^2$
tends to zero as $k$ tends to infinity.  This gives the desired inequality and
completes the proof.
\end{proof}

\bigskip
\noindent
{\bf Acknowledgements:} The author would like to thank his thesis advisor, D.H.
Phong for suggesting the problem, and for his constant advice and encouragement. 
The author is also indebted to Y-T. Siu, who visited Columbia University in the
Fall of 2002 and gave a series of lectures \cite{S5} on multiplier ideal
sheaves.  The ideas presented in those talks and in private discussions 
were invaluable in the shaping of this paper. In
particular, Professor Siu proposed the use of a complex Frobenius
theorem.    The author is also grateful to Jacob Sturm
for some very helpful suggestions, and to Zhiqin Lu for some useful discussions. 
This paper will form part of the author's forthcoming PhD thesis at Columbia
University.

\pagebreak[3]

\end{document}